\documentclass[12pt]{article}
\usepackage{amsmath,amssymb,amsthm}

\numberwithin{equation}{section}

\newtheorem{thm}{Theorem}[section]

\newcommand{\n}{\nonumber}
\renewcommand{\o}{\omega}

\renewcommand{\a}{\alpha}
\newcommand{\bb}{\begin{equation}}
\newcommand{\ee}{\end{equation}}
\newcommand{\bq}{\begin{eqnarray}}
\newcommand{\eq}{\end{eqnarray}}
\newcommand{\bqn}{\begin{eqnarray*}}
\newcommand{\eqn}{\end{eqnarray*}}

\begin{document}
\title{The global regularity for the 3D continuously stratified inviscid quasi-geostrophic equations}
\author{Dongho Chae\\
\ \\
Department of Mathematics\\
Chung-Ang University\\
 Seoul 156-756, Republic of Korea\\
email: dchae@cau.ac.kr}
\date{(to appear in {\em J. Nonlinear Science})}
\maketitle
\begin{abstract}
We prove the global well-posedness of the  continuously stratified inviscid quasi-geostrophic equations in $\Bbb R^3$.
\noindent{\bf keywords:} stratified quasi-geostrophic equations, global regularity
\end{abstract}

\section{Introduction}
\setcounter{equation}{0}
Let us consider the continuously stratified quasi-geostrophic equation for the stream function $\psi=\psi(x,y,z,t)$ on $\Bbb R^3$.
\bq\label{qg}
q_t+J(\psi, q)+\beta \psi_x =\nu\Delta q + \mathcal{F} \\
\mbox{with}\quad  q:=\psi_{xx}+\psi_{yy} +F^2\psi_{zz}.\n
\eq
Here, $F=L/L_R$ with $L$ the characteristic horizontal length of the flow and $L_R= \sqrt{gH_0}/f_0$ the Rossby deformation radius, $H_0$ the typical depth of the fluid layer and $f_0$ the rotation rate of the fluid. On the other hand, $\nu$ is the viscosity, $\mathcal{F}$ is the external force, {\em which will be set to zero for simplicity}.
In the above we used the notation, $J(f,g)=f_xg_y-f_yg_x$.  The equation (\ref{qg}) is one of the basic equations in the geophysical fluid flows.  For a physical meaning of it we mention that it can be derived from the Boussinesq equations(see\cite{ped, maj00}). In Section 1.6 of \cite{maj} one can also see a very nice explanation of (\ref{qg}) in relation to the other models of the geophysical flows.   {\em Below we consider the inviscid case $\nu=0$, and set $\beta=1$ for convenience.} The case $\nu >0$ is much easier to prove the global regularity.
Below we introduce the notations
$$ v:=v(x,y,z,t)=(-\psi_y, \psi_x, \psi_z), \quad \bar{v}:=(-\psi_y, \psi_x, 0).$$
 Rescaling in the $z$ variable as $z\to F^{-1}z$,  we have $q=\Delta \psi$. Then the equation (\ref{qg}) in our case can be written as a Cauchy problem,
\bb\label{qg1}
\left\{ \aligned &q_t +(v \cdot\nabla ) q =-v_2, \\
 & q=\Delta \psi, \\
 & v|_{t=0}=v_0.
 \endaligned \right.
\ee
Comparing the system with the vorticity formulation of the 2D Euler equations, 
\bb\label{euler} \o_t +(\bar{v} \cdot\nabla ) \o =0, \quad \o=-(\partial_x^2 +\partial_y ^2 ) \psi, \quad v=(-\psi_y, \psi_x),\ee
One can observe a similarity with the correspondence $ q \leftrightarrow \o$. We note, however,  that there exists an extra term, $v_2$ coupled, in the  first equation of (\ref{qg1}). Furthermore, more seriously,  the relation between $\psi$  and $q$ is given by a full 3D Poisson equation in the second equation of (\ref{qg1}), while in (\ref{euler}) the corresponding one is a 2D equation.
As far as the author knows the only mathematical result on the Cauchy problem of (\ref{qg1}) is 
the local in time well-posedness due to Bennett and Kloeden(\cite{ben}). 
In the viscous case there is a study of the long time behavior of solutions  of (\ref{qg}) by S. Wang(\cite{wan}).
Actually the authors 
of \cite{ben} considered 3D periodic domain  for the result, but since their proof used Kato's particle trajectory method(\cite{kat}), it is straightforward to extend the result to the case of whole domain in $\Bbb R^3$(see \cite{maj0} in the case of 3D Euler equations on $\Bbb R^3$).
In this paper our aim is to prove the following global regularity of solution to (\ref{qg1}) for  a given smooth initial data.
\begin{thm}
Let $m>7/2$, and $v_0 \in H^m (\Bbb R^3)$ be given. Then, for any $T>0$ there exists unique solution  $v\in C([0, T); H^m (\Bbb R^3))$ to the equation (\ref{qg1}).
\end{thm}
\section{Proof of the main theorem}
\setcounter{equation}{0}
\noindent{\bf Proof of Theorem 1.1 }  The local well posedness of (\ref{qg1}) for smooth $v_0$ is proved in \cite{ben}, and therefore it suffices to prove the global in time {\em a priori} estimate. Namely, we will show that
\bb\label{apriori}
\sup_{0<t<T} \|v(t)\|_{H^m}\leq C(\|v_0\|_{H^m} , T)<\infty \quad \forall T>0
\ee
for all $m>7/2$.
Taking $L^2$ inner product (\ref{qg1}) by $\psi$, and integrating by part, we obtain immediately
\bb\label{vL2}
\|v(t)\|_{L^2}=\|v_0\|_{L^2}, \quad t>0 .
\ee
Similarly,  taking $L^2$ inner product (\ref{qg1}) by $q=\Delta \psi $, and integrating by part, we obtain immediately
\bb\label{wL2}
\|q(t)\|_{L^2}=\|q_0\|_{L^2}.
\ee
Multiplying (\ref{qg1}) by $q |q|^4$, and integrating over $\Bbb R^3$, we obtain after integration by part
\bq\label{L6}
\frac{1}{6}\frac{d}{dt} \|q(t)\|_{L^6}^6 & \leq &\int_{\Bbb R^3} |v_2||q|^5 dx\leq \|v_2\|_{L^6}\|q\|_{L^6}^5\n \\
&\leq& C\|\nabla v_2\|_{L^2}\|q\|_{L^6}^5\leq C\|q\|_{L^2}\|q\|_{L^6}^5\n \\
&=& C \|q_0\|_{L^2} \|q\|_{L^6}^5,
\eq
where we used the Sobolev inequality, $\|f\|_{L^6}\leq C \|\nabla f\|_{L^2}$ in $\Bbb R^3$, and the Calderon-Zygmund estimate(see \cite{ste}), 
\bb\label{cz}\|\nabla v\|_{L^p}\leq C_p\|q\|_{L^p}, \quad 1<p<\infty.
\ee
From (\ref{L6}) we obtain
\bb
\|q(t)\|_{L^6} \leq \|q_0\|_{L^6}+ Ct \|q_0\|_{L^2}.
\ee
Hence, by the Gagliardo-Nirenberg inequality and the Calderon-Zygmund inequality we have
\bq \|v\|_{L^\infty} &\leq& C \|Dv \|_{L^6}^{\frac34} \|v\|_{L^2} ^{\frac14}\leq C \|q \|_{L^6}^{\frac34} \|v\|_{L^2} ^{\frac14}
\leq C \|q\|_{L^6} +C\|v\|_{L^2}\n \\
&\leq& C (\|q_0\|_{L^6}   +t \|q_0\|_{L^2} ) +C\|v_0\|_{L^2}  .
\eq
We introduce the particle trajectory $\{ X(a,t) \}$ on the plane generated by $\tilde{v}:=(v_1, v_2) $,
$$\frac{\partial X(a,t)}{\partial t}= \tilde{v}(X(a,t),t), \quad X(a,0)=a\in \Bbb R^2,
$$
We write (\ref{qg1}) in the form
$$
\frac{\partial}{\partial t} q(X(a,t),z,t)= -v_2 (X(a,t),z,t),
$$
which can be integrated in time as
$$
q(X(a,t),z,t) =q_0 (a, z) -\int_0 ^t v_2 (X(a,s),z,s) ds.
$$
Thus, we have
\bb\label{wLinf}
\|q(t)\|_{L^\infty} \leq \|q_0\|_{L^\infty} +\int_0 ^t \|v(s)\|_{L^\infty} ds
\leq C(1+ t^2),
\ee
 where $C=C(\|q_0\|_{L^6}, \|q_0\|_{L^2}, \|v_0\|_{L^2} )$. Combining (\ref{wL2}) and (\ref{wLinf}), using the standard $L^p$ interpolation, one has
\bb\label{wLp}
\|q(t)\|_{L^p}\leq C_1(1+t^2),\quad \forall p\in [2, \infty]
\ee
where $C_1=C_1(\|q_0\|_{L^6}, \|q_0\|_{L^2}, \|v_0\|_{L^2} )$.
 Taking $D=(\partial_1, \partial_2, \partial_3) $ on (\ref{qg1}), one has
$$
Dq_t +(D\bar{v}\cdot\nabla ) q +(\bar{v}\cdot\nabla )D q =-Dv_2.
$$
Let $p\geq 2$. Multiplying this  equation by 
$D q |Dq|^{p-2}$ and integrating it over $\Bbb R^3$, we have after integration by part, and using the H\"{o}lder
inequality and (\ref{cz}),
\bq\label{wlp1}
\frac{1}{p}\frac{d}{dt}\|Dq (t)\|_{L^p} ^p &\leq&  \| D v\|_{L^\infty} \|Dq\|_{L^p}^p +\|Dv_2\|_{L^p} \|Dq\|_{L^p}^{p-1}\n \\
&\leq& \| D v\|_{L^\infty} \|Dq\|_{L^p}^p +C \|q\|_{L^p} \|Dq\|_{L^p}^{p-1},
 \eq
from which we obtain, for $p>3$,
\bq\label{log}
\frac{d}{dt}\|Dq (t)\|_{L^p} &\leq& \| Dv\|_{L^\infty}\|Dq\|_{L^p} +C\|q\|_{L^p}\n \\
&\leq& C\{ 1+ \|Dv\|_{BMO} \log (e+\|D^2 v \|_{L^p} )\} \|Dq\|_{L^p} +C\|q\|_{L^p}\n \\
&\leq & C\{ 1+ \|q\|_{BMO} \log (e+\|Dq \|_{L^p} )\} \|Dq\|_{L^p} +C\|q\|_{L^p}\n \\
&\leq & C \{ 1+\|q\|_{L^\infty} \log (e+\|Dq \|_{L^p} )\} \|Dq\|_{L^p} +C\|q\|_{L^p}\n \\
&\leq & C (1+\|q\|_{L^\infty} + \|q\|_{L^p})( e+ \|Dq\|_{L^p} )\log (e+\|Dq \|_{L^p} ),\n \\
\eq
where we used 
the logarithmic Sobolev inequality,
\bb\label{kt}
\|f\|_{L^\infty} \leq C\{ 1+ \|f\|_{BMO} \log (e+\|Df \|_{W^{k,p}} )\}, \quad kp >3
\ee
proved in \cite{koz}, and the Calderon-Zygmund inequality.  By Gronwall's inequality,
we obtain from (\ref{log}) that
\bb\label{log1}
e+ \|Dq(t)\|_{L^p} \leq (e+ \|Dq_0\|_{L^p})^{\exp \left\{C\int_0 ^t (1+\|q(s)\|_{L^\infty} + \|q(s)\|_{L^p})ds\right\}}.
\ee
Taking into account (\ref{wLp}) and (\ref{wLinf}), we find from  (\ref{log1}) that
\bb
\sup_{0<t<T}  \|Dq(t)\|_{L^p}\leq C(v_0, T)<\infty \quad \forall T>0, \quad \forall p\in (3, \infty).
\ee
Combining  this  with the Gagliardo-Nirenberg  inequality  and (\ref{cz}), we obtain
\bq\label{gradv}
\sup_{0<t<T}  \|D v\|_{L^\infty}&\leq& C \sup_{0<t<T}  \|D^2 v\|_{L^4}  \leq C\sup_{0<t<T}  \|Dq\|_{L^4}  \n \\
&\leq &C(\|v_0\|_{W^{2,4}},  T)<\infty \quad \forall T>0.
\eq
For $p\in [2, 3]$, one has from (\ref{wlp1})  that
\bqn
\frac{d}{dt}\|Dq (t)\|_{L^p} &\leq&   \| Dv\|_{L^\infty} \|Dq\|_{L^p} +C\|q\|_{L^p} \n \\
&\leq& (\|Dv \|_{L^\infty} +\|q\|_{L^p}+1)( \|Dq \|_{L^p} +1),
\eqn
from which we obtain
\bb
\label{wlp}
\|Dq (t)\|_{L^p} +1\leq (\|Dq_0 \|_{L^p} +1)\exp \left\{ C\int_0 ^t (\|Dv (s)\|_{L^\infty} +\|q(s)\|_{L^p}+1)ds\right\}.
\ee
Hence, the estimates (\ref{wLp}) and (\ref{gradv}) imply that
\bb\label{gradw}
\sup_{0<t<T} \|Dq (t)\|_{L^p}\leq C(\|v_0\|_{W^{2,p}}, T) <\infty \quad\forall T>0.
\ee
Taking $D^2 $ on (\ref{qg}), we have
$$
D^2q_t +(D^2\bar{v}\cdot\nabla ) q +2(D\bar{v}\cdot\nabla )D q +(v\cdot\nabla ) D^2 q =-D^2v_2.
$$
Multiplying this by 
$D^2 q |D^2q|$ and integrating it over $\Bbb R^3$, we have after integration by part, and using the H\"{o}lder
inequality,
\bq\label{D2L3}
\frac{1}{3}\frac{d}{dt}\|D^2 q (t)\|_{L^3} ^3 &\leq&\|D^2 v \cdot Dq\|_{L^3} \|D^2 q\|_{L^3}^{2}
  +2 \|Dv \|_{L^\infty} \|D^2q\|_{L^3} ^3 \n \\
  &&\qquad+ \|D^2 v_2\|_{L^3} \|D^2 q\|_{L^3}^{2}\n \\
  &\leq& C(\|D^2 v\|_{BMO} \|Dq\|_{L^3} +\|D^2 v\|_{L^3} \|Dq\|_{BMO}) \|D^2 q\|_{L^3}^{2}\n \\
  &&\qquad  +2 \|Dv \|_{L^\infty} \|D^2q\|_{L^3} ^3 + \|D^2 v_2\|_{L^3} \|D^2 q\|_{L^3}^{2}\n \\
  &\leq&C\|Dq\|_{BMO}  \|D q\|_{L^3}\|D^2q\|_{L^3} ^2 \n \\
  &&\qquad  +2 \|Dv \|_{L^\infty} \|D^2q\|_{L^3} ^3 + \|Dq\|_{L^3} \|D^2 q\|_{L^3}^{2}\n \\
  &\leq&C  \|D q\|_{L^3}\|D^2q\|_{L^3} ^3 +2 \|Dv \|_{L^\infty} \|D^2q\|_{L^3} ^3\n \\
  &&\qquad   + \|Dq\|_{L^3} \|D^2 q\|_{L^3}^{2},
\eq
where we used the following bilinear estimate, proved in \cite{koz1},
$$
\|f\cdot g\|_{L^p}\leq C_p (\|f\|_{L^p} \|g\|_{BMO} +\|f\|_{BMO}\|g\|_{L^p}), \quad p\in (1, \infty)
$$
and also the critical Sobolev inequality in $\Bbb R^3$,
$$ \|f\|_{BMO}\leq C\|\nabla f\|_{L^3} .
$$
From (\ref{D2L3}) one has
\bq
\frac{d}{dt}\|D^2q (t)\|_{L^3} &\leq& C ( \|D q\|_{L^3} +  \|Dv \|_{L^\infty} ) \|D^2q\|_{L^3}+\|Dq\|_{L^3} \n \\
 &\leq &  C ( \|D q\|_{L^3} +  \|Dv \|_{L^\infty} +1) (\|D^2q\|_{L^3} +1),
 \eq
 from which one has
 \bb
 \|D^2 q(t)\|_{L^3}+1 \leq (\|D^2 q_0\|_{L^3} +1)\exp \left\{ C\int_0 ^t  ( \|D q(s)\|_{L^3} +  \|Dv (s)\|_{L^\infty} +1)ds \right\}.
 \ee
 Combining this with (\ref{gradw}) and (\ref{gradv}), we have
 \bb\label{bmo}
 \sup_{0<t<T}\|Dq(t)\|_{BMO}\leq C\sup_{0<t<T} \|D^2 q(t)\|_{L^3} \leq C(\|v_0\|_{W^{3,3}}, T) <\infty \quad \forall T>0.
 \ee
 Let $\a=(\a_1, \a_2,\a_3)\in (\Bbb N\cup \{0\})^3$ be a muti-index. Let $m>7/2$.
 Taking $D^\a $ on (\ref{qg1}), and multiplying it by $D^\a q$,
 summing over $|\a|\leq m-1$, and integrating it over $\Bbb R^3$, we have
 \bq\label{gen}
 \frac12 \frac{d}{dt} \|q(t)\|_{H^{m-1}} ^2&=& \sum_{|\a|
 \leq m-1} (D^\a (\bar{v}\cdot \nabla )q , D^\a q )_{L^2}  -\sum_{|\a|\leq m-1} (D^\a v_2, D^\a  q)_{L^2} \n \\
 &=& \sum_{|\a|
 \leq m} (D^\a (\bar{v}\cdot \nabla )q-(\bar{v}\cdot \nabla ) D^\a q , D^\a q )_{L^2}  -\sum_{|\a|\leq m-1} (D^\a v_2, D^\a  q)_{L^2} \n \\ 
 &\leq&\sum_{|\a|\leq m-1} \|D^\a (\bar{v}\cdot \nabla )q-(\bar{v}\cdot \nabla ) D^\a q \|_{L^2}\|q\|_{H^{m-1}} + \|v\|_{H^{m-1}} \|q\|_{H^{m-1}}\n \\
 &\leq& C(\|\nabla v\|_{L^\infty}+ \|\nabla q\|_{L^\infty}) (\|q\|_{H^{m-1} } + \|v\|_{H^{m-1}})\|q\|_{H^{m-1}}+ \|v\|_{H^{m-1}} \|q\|_{H^{m-1}},\n
  \\
 \eq
 where we used the following commutator estimate,
 $$
 \sum_{|\a|\leq m} \|D^\a (fg)-fD^\a g\|_{L^2}\leq C_m (\|\nabla f\|_{L^\infty} \|D^{m-1} g\|_{L^2} 
 +\|D^m f \|_{L^2} \|g\|_{L^\infty} ),
 $$
 proved in \cite{kla}.
 We observe the following norm equivalence: there exists  a constant $K$  independent of $v,q$  such that 
 $$ K^{-1} (\|v\|_{L^2}^2+ \|q\|_{H^{m-1}} ^2)\leq \|v\|_{H^m} ^2\leq K(\|v\|_{L^2}^2+ \|q\|_{H^{m-1}}^2 ),
 $$
 which is an immediate consequence of (\ref{cz}).
  Since $\|v(t)\|_{L^2}=\|v_0 \|_{L^2}$ as in (\ref{vL2}), one can add $\|v_0\|_{L^2}^2$  to $\|q\|_{H^{m-1}}^2$ in the left hand side of (\ref{gen}) to obtain
  \bq
\frac{d}{dt}  Y(t)&\leq& C (\|\nabla v\|_{L^\infty} +\|\nabla q\|_{L^\infty}+1 )Y(t)\n \\
 &\leq&  C(\|\nabla v\|_{BMO} +\|\nabla q\|_{BMO} +1)Y(t) \log Y(t) \eq
where we set
$$
Y(t):=e+\|v(t)\|_{L^2}^2 +\|q(t)\|_{H^m}^2,
$$
and used  the logarithmic Sobolev inequality in the form (\ref{kt}) for $m>7/2$.   
By Gronwall's inequality we have
 \bq\label{last}
 e+\|v(t)\|_{H^m} ^2  &\leq& e+C(\|v(t)\|_{L^2}^2+ \|q(t)|_{H^{m-1}} ^2) \n \\
 &\leq& e+C(\|v_0\|_{L^2}^2+ \|q_0\|_{H^{m-1}} ^2) ^{\exp\left\{ C\int_0 ^t (\|\nabla v(s)\|_{BMO} +\|\nabla q(s)\|_{BMO} +1) ds\right\}}\n \\
 &\leq&e+ (C\|v_0\|_{H^m} ^2 )^{\exp\left\{ C\int_0 ^t (\|\nabla v(s)\|_{L^\infty} +\|\nabla q(s)\|_{BMO} +1) ds\right\}}.\n \\
 \eq
 The estimates (\ref{gradv}) and (\ref{bmo}), combined with (\ref{last}), provides us with (\ref{apriori}). $\square$

 $$ \mbox{\bf Acknowledgements } $$
 The author would like to thank deeply to Prof. Shouhong Wang for suggesting the problem
 and helpful discussions.
This research is supported partially by NRF
  Grants no.2006-0093854 and  no.2009-0083521.

\end{document}